\documentclass[11pt, reqno]{amsart}
\usepackage{diagbox, amsmath, amsthm, amscd, amsfonts, amssymb, graphicx, color, pgf, tikz}
\usepackage{booktabs}
\usepackage[bookmarksnumbered, colorlinks, plainpages]{hyperref}
\usepackage[all,cmtip]{xy} 
\usepackage{mathrsfs}
\usetikzlibrary{arrows}
\usepackage{amsmath}
\usepackage[autostyle]{csquotes}

\makeatletter \oddsidemargin.9375in \evensidemargin \oddsidemargin
\newtheorem{theorem}{Theorem}[section]
\newtheorem{lemma}[theorem]{Lemma}
\newtheorem{proposition}[theorem]{Proposition}
\newtheorem{corollary}[theorem]{Corollary}
\theoremstyle{definition}
\newtheorem{definition}[theorem]{Definition}
\newtheorem{example}[theorem]{Example}
\newtheorem{question}[theorem]{Question}

\newtheorem{remark}[theorem]{Remark}
\numberwithin{equation}{section}

\begin{document}
\title[Countable-configurations and paradoxical decompositions]{Countable-configurations and countable paradoxical decompositions}

\author{M. Meisami, A. Rejali$^{*}$ and A.Yousofzadeh}

\address{ Mahdi Meisami, Department of Pure Mathematics, Faculty of Mathematics and Statistics, Ph.D student, University Of Isfahan , Isfahan 81746-73441, IRAN.}\email{m.meisami@sci.ui.ac.ir}
\address{Ali Rejali, Department of Pure Mathematics, Faculty of Mathematics and Statistics, University Of Isfahan, Isfahan 81746-73441, Iran}\email{rejali@sci.ui.ac.ir , alirejali04@gmail.com}

\address{Akram Yousofzadeh, Department of Mathematics, Mobarakeh Branch Islamic Azad University, Isfahan, Iran}\email{ayousofzade@yahoo.com}
\begin{abstract}
In this paper we define countable-configuration of groups and prove that two Hopfian groups with the same set of countable-configurations are isomorphic and vice versa. 
We also study the countable paradoxical decomposition of groups. It is proved that a group $G$ admits a countable paradoxical decomposition if and only if it is infinite. 
\end{abstract}

\keywords{Amenability, paradoxical decomposition, group action, Configuration.\\
\indent $^{*}$ Corresponding author}

\maketitle


\newcommand\sfrac[2]{{#1/#2}}

\newcommand\cont{\operatorname{cont}}
\newcommand\diff{\operatorname{diff}}

\section{\bf Introduction}

The concept of configuration for groups was introduced by Rosenblatt and Willis. They used this concept particularly to characterize the amenability of groups \cite{rw}. Configurations are also applied to construct paradoxical decompositions and to estimate Tarski numbers of non-amenable groups \cite{ry2,yyy,yme,rty}. This notion also was generalized to the case of semigroups \cite{ar}.

The Banach-Tarski theorem states that a solid 3-dimensional ball can be decomposed into finitely many pieces in such a way that after reassembling through translations and rotations only, two identical balls can be obtained. The essence of the proof is that the rotation group $SO_3(\mathbb{R})$ contains a copy of the free group $\mathbb{F}_2$ with two generators. As $\mathbb{F}_2$ admits a free subgroup of denumerable rank, one can give a similar partition of the unit sphere in the 3-dimensional Euclidean space with infinitely countably many pieces.

It should be noted that the Banach-Tarski paradox is valid even continuously. Churkin shows this fact by proving that $SO_3(\mathbb{R})$ contains a free group whose rank is the power of continuum (see \cite{chur},\cite{tom}).

The main subject of this paper is to extend configurations and paradoxical decompositions of groups and group actions with emphasis on countable subsets. So recalling the notion of finite configurations and finite paradoxical decompositions turns out to be necessary. We first mention the definition and notations we used throughout the paper.  For the motivation behind configuration see \cite{sur}.

If $G$ is a group with identity element $e_G$, and $X$ is a set, then a left group action of $G$ on $X$ is a function
 $\cdot :G\times X\longrightarrow X$ such that $e_G\cdot x=x$ , and $g_1\cdot(g_2\cdot x)=(g_1g_2)\cdot x$ , for all $g_1,g_2\in G$ and $x\in X$. We often denote a group action of $G$ on $X$ by $G\curvearrowright X$ if the action is fixed.
   The orbit of $x\in X$ which is denoted by $[x]$ is defined as $[x]=\{g.x \mid\ \ g\in G\}.$ 
Clearly the set of orbits of elements of $X$ forms a partition of X. 
We say that $G$ acts transitively on $X$ if for every $x,y\in X$ there is an element $g\in G$ such that $g\cdot x=y$, so the action admits only one orbit. The notation $gx$ is used instead of $g.x$ when the action being considered is clear from context.
The set of all orbits of $X$ under the action of $G$ is written as $X/G$ and is called the quotient of the action $G\curvearrowright X$.
The stabilizer  of $x\in X$ which is a subgroup of $G$ is the set of all elements in $G$  fixing $x$,
$$G_{x}=\{g\in G\mid g\cdot x=x\}.$$

Let $G\curvearrowright X$ be a group action, $\mathcal E=\{E_1,\dots, E_m\}$ be a finite partition of $X$ and $\mathfrak g=(g_1,\dots,g_n)$ be a finite ordered subset of $G$.  An $(n+1)$-tuple 
$C=(C_0,C_1,\dots,C_n)$ of integer numbers in $\{1,\dots, m\}$ is called a configuration corresponding to the pair $(\mathfrak g,\mathcal E)$ if there exists $x_0\in X$ such that for 
each $j\in\{1,\dots,n\}$ we have $g_jx_0\in E_{C_j}$. The set of all such configurations is denoted by $Con(\frak g, \mathcal E ;X)$ and the family of all possible  $Con(\frak g, \mathcal E;X)$ is denoted by $Con(G\curvearrowright X)$.  Two group actions $G\curvearrowright X$ and $H\curvearrowright Y$ are called configuration equivalent if $Con(G\curvearrowright X)=Con(H\curvearrowright Y)$. We use the notations $Con(\frak g, \mathcal E)$ and $Con(G)$ when $G$ acts on itself by left multiplication.

If in above definition, the ordered set $\mathfrak g$ is assumed to be a generating set for $G$, then we use the notations $Con^*(\frak g, \mathcal E)$ and $Con^*(G)$ in place of $Con(\frak g, \mathcal E)$ and $Con(G)$, respectively. We say two groups $G$ and $H$ are $*$-configuration equivalent if $Con^*(G)=Con^*(H)$. Such groups share various group-theoretical properties. For more details see \cite{arw,ary,sr1,sr2,ry2}.
It is observed that these two notions are different. For example if $H$ is any arbitrary subgroup of $\mathbb F_n$ (the non-amenable free group on $n$ generators), then $Con(H)\subseteq Con(\mathbb F_n)$. But $Con^*(H)\subseteq Con^*(\mathbb F_n)$ only if $H\cong \mathbb F_n$ (see \cite{arw} and \cite{mrsy}).

Let $G$ be a discrete group, $\{A_1,\dots,A_n, B_1,\dots,B_m\}$ be a finite partition for $G$ and $\{g_1,\dots,g_n, h_1,\dots, h_m\}$ be a subset of $G$. Then $(A_i, B_j, g_i, h_j)_{i,j=1}^{n,m}$ is called a paradoxical decomposition for $G$ if 
$$G=\left (\bigsqcup_{i=1}^nA_i\right)\bigsqcup\left (\bigsqcup_{j=1}^m B_j\right)=\bigsqcup_{i=1}^ng_iA_i=\bigsqcup_{j=1}^mh_jB_j.$$

 By Tarski alternative a discrete group is amenable if and only if it admits no  (finite) paradoxical decompositions. If the cardinal number of pieces in the definition of a paradoxical decomposition is allowed to be $\aleph_0$, one faces the concept of countable (or geometric) paradoxical decomposition. There are interesting facts about this decomposition. For instance, the amenable group $(\mathbb R,+)$ admits a countable paradoxical decomposition. For some detailed results about this kind of decompositions see  \cite[P.7]{wagon}.  

The initial point in the theory of paradoxical decomposition may be returned to Galileo, who observed that the set $\mathbb{N}$ and its subset of squares have the same cardinality. After that set theorists have found out that each infinite set $X$ satisfies $\vert X\vert=2\vert X\vert$. Then the theory of Banach-Tarski began and different types of finite paradoxical decompositions were introduced. In this paper, we are interested to investigate various properties of countable paradoxical decompositions.\\

In section 2 we introduce countable-configurations for groups and group actions and prove that two countable-configuration equivalent groups are isomorphic under a condition.
In section 3 we devote ourselves to countable paradoxical decompositions of  group actions. 


\section{countable-configuration}
Suppose that $G\curvearrowright X$ is a group action, $G$ is a countable group and $X$ is an infinite set. Let  $\mathfrak{g}=\{g_i\}_{i=1}^\infty$ be a sequence of elements of $G$ and $\mathcal{E}=\{E_i\}_{i=1}^\infty$ be an infinite partition for $X$. By a countable-configuration of this group action related to $(\mathfrak{g},\mathcal{E})$ we mean a sequence $C=(C_0,C_1,\ldots)$ of positive integers such that  for some $x\in E_{C_0}$, 
$$g_i \cdot x \in E_{C_i},\quad  (i \in \mathbb{N}).$$ 
We say that $(\mathfrak{g},\mathcal{E})$ is a countable-configuration pair if $\mathfrak{g}$ is an infinite sequence in $G$ and $\mathcal{E}$ is an infinite partition of $X$. 
We will denote the set of all countable-configurations of group action related to $(\mathfrak{g},\mathcal{E})$ by $\mathrm{Con}_c(\mathfrak{g},\mathcal{E};X)$. Also we define 
\begin{equation*}
\mathrm{Con}_c(G\curvearrowright X)=\{\mathrm{Con}_c(\mathfrak{g},\mathcal{E};X)\mid (\mathfrak{g},\mathcal{E})\;\text{is a countable-configuration pair}\}.
\end{equation*}

For a  configuration  $C\in Con_c(\mathfrak{g},\mathcal{E};X)$ we set 
$$x_0(C)=E_{C_0}\cap(\cap_{i=1}^\infty g_i^{-1}E_{C_i})$$
and 
$$x_j(C)=g_jx_0(C)\quad \  (j\geq 1).$$
According to the above definition, we have
$$C=(C_0,C_1,\ldots)\in \mathrm{Con}_c(\mathfrak{g},\mathcal{E};X)\text{ if\ and\ only\ if }x_0(C)\neq \emptyset.$$
Observe that 
\begin{equation}\label{efraz} 
E_i=\bigsqcup_{C_j=i}x_j(C),\qquad 1\leq i\leq m.
\end{equation}
Then we can easily see that for $j\geq 0$, the familiy $\{x_j(C): \ C\in Con_c(\frak g, \mathcal E)\}$ makes an infinite partition for $X$.

As in the definition of $Con_c(\frak g, \mathcal E;X)$ it is emphasized that  $G$ and $X$ are infinite sets, we need to fix the convention $\mathrm{Con}_c(G\curvearrowright X)=\mathrm{Con}(G\curvearrowright X)$, when $G$ or $X$ is finite.
%
\begin{proposition}\label{imp}
Let $G$ and $H$ be two countable groups acting on the sets $X$ and $Y$ respectively and $Con_c(G\curvearrowright X)=Con_c(H\curvearrowright Y).$ Then 
$$Con(G\curvearrowright X)=Con(H\curvearrowright Y).$$ 
\end{proposition}
\begin{proof}
It is trivial, if one of the sets $G$ or $X$ in the statement is finite. Let $G$ and $X$ be infinite and $Con_c(G\curvearrowright X)=Con_c(H\curvearrowright Y).$ Then clearly $Y$ is infinite too. 
 Suppose that $Con(\mathfrak{g},\mathcal{E};X)\in Con(G\curvearrowright X)$ where $\mathfrak{g}=(g_1,\ldots,g_n)$ and $\mathcal{E}=\{E_1,\ldots,E_m\}$.  
 Since $X$ is infinite, we can assume without loss of generality that $E_m$ is infinite. Let $E_m=\{y_m,y_{m+1},\dots\}$  
 and $E'_j=\{y_j\}$, for $j\geq m$. Define 
\begin{equation*}
{\mathfrak{g'}}=(g_1,\ldots , g_n,e,e,\ldots ) \text{ and } {\mathcal{E'}}=\{E_1,\ldots,E_{m-1}, E'_{m}, E'_{m+1},\dots\}.
\end{equation*}
Then
$$
Con({\mathfrak{g'}},{\mathcal{E'}};X)\in Con_c(G\curvearrowright X)=Con_c(H\curvearrowright Y).
$$
So there is a configuration pair $({\mathfrak{h'}},{\mathcal{F'}})$ such that $$Con({\mathfrak{g'}},{\mathcal{E'}};X)=Con({\mathfrak{h'}},{\mathcal{F'}};Y),$$ for some string ${\mathfrak{h'}}=(h_1,h_2,\ldots )$ of elements of $H$ and partition ${\mathcal{F'}}=\{F'_1,F'_2,\ldots \}$ of $Y$. \\
Note that  $Con({\mathfrak{g'}},{\mathcal{E'}};X)$ consists of configurations of the form
$D=(d_0,d_1,\dots)$ for which $d_j=d_0$ if $j\geq m$ and there exists $C=(C_0,C_1,\dots,C_n)\in Con(\frak g,\mathcal E;X)$ such that $c_k=m$ if and only if $d_k\geq m, \quad (0\leq k\leq n) .$
Therefore      $Con({\mathfrak{g}},{\mathcal{E}};X)=Con({\mathfrak{h}},{\mathcal{F}};Y)$, where $\frak h=(h_1,\dots,h_n)$ and $\mathcal F$ is the partition $\{F'_1,\dots,F'_{m-1},F_m\}$  in which $F_m=\cup_{k=m}^\infty{F'_k}.$ \\
Since $\mathfrak{g}$ and $\mathcal{E}$ were chosen arbitrary, we have  $Con(G\curvearrowright X)\subseteq Con(H\curvearrowright Y)$. Similarly $Con(H\curvearrowright Y)\subseteq Con(G\curvearrowright X)$.
\end{proof}
 Proposition \ref{imp} implies that every result that has been obtained for configuration equivalent groups, also holds for countable-configuration equivalent groups. 
 In particular such groups have the same Tarski numbers (see \cite{ry2}).

We need the following  definition for the second part of Theorem  \ref{last} below. 

\begin{definition}
The group  $G$ is called a $Q$-group (Hopfian) if every epimorphism of $G\to G$ is an isomorphism.
\end{definition}
Equivalently $G$ is a $Q$-group if for every normal subgroup  $N$ of $G$ such that $G$ and $G/N$ are isomorphic,
we have $N=\{e_G\}$. For more details about $Q$-groups, the reader is referred particularly to \cite{baer}. We now state the main result of this paper.

\begin{theorem}\label{last}
 If $G$ and $H$ are two at most countable groups with $Con^*_c(G)=Con^*_c(H)$, then $H$ admits a normal subgroup $F$ such that $G\cong \frac{H}{F}$.
 If in addition $G$ is a $Q$-group, then $G\cong H$.
\end{theorem}
\begin{proof}
If $G$ is finite, by our convention $Con(G)=Con(H)$ which implies that $G\cong H$ by \cite[Proposition 4.3]{arw}.  
Let $G$ be infinite and $Con^*_c(G)=Con^*_c(H)$. By assumption $G$ is countable, say $G=\{g_i\mid  i\in \mathbb N\}$, where $g_1=e_G$, the identity element of $G$. Set $\frak g=(g_1,g_2,\dots)$ and $\mathcal E=\{E_1,E_2,\dots\}$ for which $E_j=\{g_j\}$. Suppose that $\pi:\mathbb N\times \mathbb N\rightarrow \mathbb N$ be the function defined through $\pi(i,j)=k$ if and only if $g_ig_j=g_k.$
Then $Con(\frak g,\mathcal E)$ consists of configurations of the form
\begin{equation}\label{cj}
C_j=(j, \pi(2,j), \pi(3,j),\dots), \quad (j\in \mathbb N).
\end{equation}

Since by assumption $Con^*_c(G)=Con^*_c(H)$, there exists a configuration pair $(\frak h,\mathcal F)$ for $H$ such that 
$Con^*_c(\frak g,\mathcal E)=Con^*_c(\frak h,\mathcal F)$. Then for each $j\in \mathbb N$, $C_j\in Con^*_c(\frak h,\mathcal F).$ This implies that for each $j,\ell\in \mathbb N$ we have $h_jF_\ell\subseteq F_{\pi(j,\ell)}.$ Without loss of generality we assume that $e_H\in F_1.$
 By (\ref{cj}) $x_{j-1}(C_\ell)\subseteq F_{\pi(j,\ell)}$ 
and for every $\ell'\neq \ell,$ 
$x_{j-1}(C_{\ell'})\nsubseteq F_{\pi(j,\ell)}$ and on the other hand $\{x_{j-1}(C_\ell)\mid \ell\in \mathbb N\}$ is a partition for $H$. Hence $x_{j-1}(C_\ell)=F_{\pi(j,\ell)}$. In particular 
\begin{equation}\label{eq}
h_jF_{\ell}=x_{j-1}(C_\ell)=F_{\pi(j,\ell)}.
\end{equation}


\noindent{Now we prove that $F_{1}$ is a normal subgroup of $H$.}
\\
 
 \noindent  Let $\pi(j,j')=1.$ By (\ref{eq}), $h_jF_{j'}=F_{\pi(j,j')}=F_1$. Thus 
 $$h_j^{-1}F_{1}=F_{j'}=h_{j'}F_1.$$
Now let $e_H\neq x\in F_1.$ Since $H$ is generated by $\frak h,$ $x$ is a finite combination of elements of $\frak h\cup \frak h^{-1},$ say, 
$x=h_{i_1}^{\alpha_{i_1}}\cdot \dots \cdot h_{i_p}^{\alpha_{i_p}}$ for some $i_1,\dots i_p\in \mathbb N$
and some $\alpha_{i_1},\dots \alpha_{i_p}\in \{1,-1\}.$
 It is then readily seen that for  $t=\pi(j_1,\dots,\pi(j_{p-1},j_p))$, where
 $$j_q=\begin{cases}
 i_q\ \ \ \text{if}\ \alpha_{i_q}=1\\
 i'_{q}\ \ \ \text{if}\ \alpha_{i_q}=-1
 \end{cases}.$$
$xF_1=h_tF_1=F_t$. On the other hand $e_H\in F_1.$ So, 
$$x\in F_1\cap F_t,$$
which is possible only if $t=1.$ It means in particular that $F_1$ is a subgroup of $H$.
 To prove the normality of $F_1$ it is enough to show that the set $\{F_j\mid j\geq 1\}$ of left cosets of $F_1$ in $H$  is closed under multiplication. Let $x=h_{i_1}^{\alpha_{i_1}}\cdot \dots \cdot h_{i_p}^{\alpha_{i_p}}\in F_1$ for some $i_1,\dots i_p\in \mathbb N$. Then $h_kxh_jF_1=h_{\pi(k,\pi(t,j)})F_1=F_{\pi({k,j})}$. Therefore $F_kF_j=F_{\pi(k,j)},$ for each $k,j\in \mathbb N$. This implies the normality of $F_1$.\\
  Finally the function $\phi:h_j\mapsto F_j$ is an isomorphism from $G$ onto $\frac{H}{F_{1}}$. \\
 
 \noindent 

Now let in addition $G$ be a $Q$-group. For the second part of the theorem, it is enough to show that $F_1=\{e_H\}$. Suppose that $F_1\neq \{e_H\}$ and set $$F^1_1=\{e_H\} \ \ \ \text{and}\ \ \ \ \ F^2_1=F_1\setminus \{e_H\}.$$
Also set $$\mathcal F_1=\{F_1^1, F_1^2, F_2, F_3,\dots\}.$$
Then there exists a countable-configuration pair $(\mathfrak g_1,\mathcal E_1)$ for $G$ such that $Con^*_c(\frak g_1,\mathcal E_1)=Con^*_c(\frak h,\mathcal F_1)$.
Let $\mathcal E_1=\{E_1^1, E_1^2, E_2, E_3,\dots\}$ and $e_G\in E_1^1$. By a similar argumentation as above, ${E^1_1\cup E_1^2}$ is a normal subgroup of $G$  and $G\cong \frac{G}{E^1_1\cup E^2_2}.$ Since $G$ is assumed to be a $Q$-group, we have $E_1^2=\emptyset$ which is impossible. This implies that $F_1^2=\emptyset$ and the proof is complete.

\end{proof}
\begin{question}
Let $Con^* (G)=Con^* (H)$ and $G$ is a $Q$-group. Does $H$ is a $Q$-group and $G\cong H$?
\end{question}

\section{Countable decompositions}

 Let $G$ be a group acting on a set $X$, and consider two subsets $A$ and $B$ of $X$. The sets
$A$ and $B$ are finitely $G$-equidecomposable (resp. countably $G$-equidecomposable) if there exist  finite (resp. infinite) sequences  of pairwise disjoint subsets $\{A_i\}_{i=1}^n$ and $\{B_i\}_{i=1}^n$ (resp. $\{A_i\}_{i=1}^{\infty}$ and $\{B_i\}_{i=1}^{\infty}$) and group elements $\{g_i\}_{i=1}^n$ (resp. $\{g_i\}_{i=1}^{\infty}$) such that $A=\bigcup_{i=1}^n A_i$ and $B=\bigcup_{i=1}^n B_i$ (resp. $A=\bigcup_{i=1}^{\infty} A_i$ and $B=\bigcup_{i=1}^{\infty} B_i$) and $B_i=g_iA_i$ for all $i\in \{1,\cdots ,n\}$ (resp. $i\in \mathbb N$).

   If $A$ and $B$ are finitely $G$-equidecomposable, we write $A\sim_G B$ and if $A$ and $B$ are countably $G$-equidecomposable, we write $A\overset{c}{\sim}_G B$. If the group under consideration is definite from the context, then we can omit $G$ and write $A\sim B$ and $A\overset{c}{\sim} B$. 
 \begin{definition}
If $G$ is a group acting on a set $X$, then
\begin{itemize}
\item[(1)] A subset $A$ of $X$ is finitely $G$-paradoxical if there exist two subsets of $A$ named $M$ and $N$ such that $A=M\bigsqcup N$ and $M\sim_G A\sim_G N$.
 \item[(2)] A subset $A$ of $X$ is countably $G$-paradoxical if there exist two subsets of $A$ named $M$ and $N$ such that $A=M\bigsqcup N$ and $M\overset{c}{\sim}_G A\overset{c}{\sim}_G N$.
\end{itemize}
 \end{definition}
 When $G\curvearrowright G$ with the group operation from left and we have a subset $A$ of the group that is countably (finitely) $G$-paradoxical, simply we say that $A$ is countably (finitely) paradoxical.
 
 Note that every infinite countable group $G$ is countably paradoxical. To see this, it is enough to dispart $G$ into two disjoint infinite countable subsets and write each of them as the disjoint union of singletons.
 
\begin{lemma}\label{sub}
If the group $G$ has a countably paradoxical subgroup $H$, then $G$ is countably paradoxical too.
\end{lemma}
\begin{proof}
Since $H$ is countably paradoxical,  it admits a countable partition $\{A_i\}_{i=1}^{\infty}\sqcup \{B_i\}_{i=1}^{\infty}$ such that $H=\bigsqcup _{i=1}^\infty g_i A_i=\bigsqcup_{i=1}^\infty h_i B_i$. Choose a right transversal of $H$ in $G$, namely, a subset $A$ of $G$ such that $G=\bigsqcup_{x\in A} Hx$. Then
\begin{eqnarray*}
G&=&\bigsqcup_{x\in A} Hx=\bigsqcup_{x\in A}  \left( \left(\bigsqcup_{i=1}^\infty A_i\right)\bigcup \left(\bigsqcup_{i=1}^\infty B_i\right)\right) x\\
&=&\left( \bigsqcup_{i=1}^\infty \bigsqcup_{x\in A} A_i x\right) \bigcup \left( \bigsqcup_{i=1}^\infty \bigsqcup_{x\in A} B_i x\right)\\
&=&\left( \bigsqcup_{i=1}^\infty  A_i A\right) \bigcup \left( \bigsqcup_{i=1}^\infty  B_i A\right).
\end{eqnarray*}
According to our assumption we have analogously 
\begin{eqnarray*}
G&=&\bigsqcup_{i=1}^\infty \bigsqcup_{x\in A} (g_i A_i)x=\bigsqcup_{i=1}^\infty \bigsqcup_{x\in A} (h_i B_i)x\\
&=&\bigsqcup_{i=1}^\infty g_i A_i A=\bigsqcup_{i=1}^\infty h_i B_i A.
\end{eqnarray*}
Therefore the family $\{A_i A\}_{i=1}^\infty\cup\{B_i A\}_{i=1}^\infty$ along with $\{g_i\}_{i=1}^\infty\cup\{h_i\}_{i=1}^\infty$ forms a countable paradoxical decomposition for $G$.
\end{proof}
 \begin{theorem}\label{finite}
A group $G$ is countably paradoxical if and only if it is infinite. 
\end{theorem}
\begin{proof}
 Clearly no finite group admits a countable paradoxical decomposition. On the other hand every infinite group possesses an infinite countable subgroup (for example the subgroup generated by a countable subset) which is countably paradoxical. Now Lemma \ref{sub} implies that $G$ is countably paradoxical too.
\end{proof}

\begin{example}
Let $G\curvearrowright X$ be the trivial action, i.e. $g\cdot x=x$, for each $g\in G$ and $x\in X$. This action is not countably paradoxical. Suppose it is not true. Then $G\curvearrowright X$ admits a countable paradoxical decomposition. In other words there exist an infinite partition $\{A_i \}_{i=1}^\infty\cup\{B_i \}_{i=1}^\infty$ for $G$ and   $\{g_i\}_{i=1}^\infty\cup\{h_i\}_{i=1}^\infty\subseteq G$ such that
$$
X=\bigsqcup_{i=1}^\infty g_i A_i=\bigsqcup_{i=1}^\infty A_i
$$
Similarly
$$
X=\bigsqcup_{j=1}^\infty h_j B_j=\bigsqcup_{j=1}^\infty B_j
$$
which contradicts 
$$
X=\left( \bigsqcup_{i=1}^\infty A_i \right) \bigsqcup \left( \bigsqcup_{j=1}^\infty B_j \right).
$$
\end{example}

\begin{corollary}
If $G$ is not countably paradoxical, then each action $G\curvearrowright X$ is not countably paradoxical as well.
\end{corollary}
\begin{proof}
Suppose that there is an action $G\curvearrowright X$ which is countably paradoxical. Then one finds a familiy $\{A_i\}_{i=1}^\infty \sqcup\{B_j \}_{j=1}^\infty$ of disjoint subsets of $X$ and also two sequences of elements $\{g_i\}_{i=1}^\infty$ and $\{h_j \}_{j=1}^\infty$ in $G$ such that
$$X=\left( \bigsqcup_{i=1}^\infty A_i \right) \bigsqcup \left( \bigsqcup_{j=1}^\infty B_j \right)=\bigsqcup_{i=1}^\infty g_i A_i = \bigsqcup_{j=1}^\infty h_j B_j.$$\\
If $G$ is not countably paradoxical, then it is finite by Theorem \ref{finite}. So there are  $n,m\in \mathbb N$ and $$g_1',\dots,g_n',h_1',\dots,h_m'\in G$$ such that  $\{g_i \}_{i=1}^\infty=\{g_1',\dots,g_n'\}$ and $\{h_i \}_{i=1}^\infty=\{h_1',\dots,h_m'\}$. Then we can write the decomposition of $X$ as follow 
$$
X=\bigsqcup_{i=1}^\infty g_i A_i =\bigsqcup_{j=1}^{n} g_j' A_j' $$
and
$$
X=\bigsqcup_{i=1}^\infty h_i B_i =\bigsqcup_{k=1}^{m} h_k' B_k', $$
where
$$ A_j':=\bigsqcup_{i=1}^\infty \{ A_i : g_i=g_j'\}$$
and
$$B_k':=\bigsqcup_{i=1}^\infty \{ B_i : h_i=h_k'\}.$$
In other words, $G\curvearrowright X$ is finitely paradoxical which is impossible since every finite group is amenable and no action of an amenable group  is paradoxical by \cite[Lemma 2.3]{golan}.
\end{proof}


If $G$ acts trivially on  an infinite countable set $X$, $\vert X/G \vert = \vert X \vert$.  So $X/G$ is an infinite countable set, as well. On the other hand we know that this action does not admit a countable paradoxical decomposition. In the following theorem we give a necessary and sufficient condition for  $G\curvearrowright X$  to be countably paradoxical, when the quotient of this action is countable.
\begin{theorem}\label{thinfinite}
Let $X/G$ be infinite and countable. Then $G\curvearrowright X$ has countable paradoxical decomposition if and only if for every $x\in X$, $G/G_x$ is ?countably $?G?$?-paradoxical.
\end{theorem}
\begin{proof}
Suppose that each $G/G_x$ is countably $?G?$?-paradoxical. So for each $x\in X$ there are subsets $\{A_i^x\}_{i=1}^\infty , \{B_j^x\}_{j=1}^\infty$ of $G/G_x$ ???and elements $\{g_i^x\}_{i=1}^\infty , \{h_j^x\}_{j=1}^\infty$ in  $G$ such that 
$$
G/G_x=\left( \bigsqcup_{i=1}^\infty A_i^x \right) \bigsqcup \left( \bigsqcup_{j=1}^\infty B_j^x \right) =  \bigsqcup_{i=1}^\infty g_i^x A_i^x = \bigsqcup_{j=1}^\infty h_j^x B_j^x .
$$
Let $\phi_x : G/G_x\longrightarrow [x]$ be such that $hG_x\longmapsto hx$. Then $\phi_x$ is a bijection and translation invariant.
Thus
$$
[x]=\left( \bigsqcup_{i=1}^\infty \phi_x (A_i^x) \right) \bigsqcup \left( \bigsqcup_{j=1}^\infty \phi_x(B_j^x) \right) =  \bigsqcup_{i=1}^\infty g_i^x \phi_x (A_i^x) = \bigsqcup_{j=1}^\infty h_j^x \phi_x (B_j^x) .
$$
Since $X=\bigsqcup_{x\in E} [x]$ for some countable set $E$,

\begin{align*}
X=& \left( \bigsqcup_{i=1}^\infty \bigsqcup_{x\in E} \phi_x (A_i^x) \right) \bigsqcup \left( \bigsqcup_{j=1}^\infty \bigsqcup_{x\in E} \phi_x(B_j^x) \right)
\\&=\bigsqcup_{i=1}^\infty \bigsqcup_{x\in E} g_i^x \phi_x (A_i^x)
=\bigsqcup_{j=1}^\infty \bigsqcup_{x\in E} h_j^x \phi_x (B_j^x) .
\end{align*}
Conversely if the action $G\curvearrowright X$ has countable paradoxical decomposition we can write 
$$
X=\left( \bigsqcup_{i=1}^\infty A_i \right) \bigsqcup \left( \bigsqcup_{j=1}^\infty B_j \right) =  \bigsqcup_{i=1}^\infty g_i A_i = \bigsqcup_{j=1}^\infty h_j B_j .
$$
Hence for every $x\in X$, by intersection, we achieve that
$$
[x]=\left( \bigsqcup_{i=1}^\infty A_i \cap [x] \right) \bigsqcup \left( \bigsqcup_{j=1}^\infty B_j \cap [x] \right) =  \bigsqcup_{i=1}^\infty (g_i A_i \cap [x])= \bigsqcup_{j=1}^\infty (h_j B_j \cap [x]) .
$$
For the second equation note that for every $g\in G$ and $x\in X$, $g[x]=[x]$. In particular
$$
g_i A_i \cap [x]=g_i A_i \cap g_i [x]=g_i (A_i\cap [x]).
$$ 
 As a result $[x]$ (as a subset of $X$) has a countable $G$-paradoxical decomposition. Since the function $\phi_x$ is bijection and translation invariant, the set $G/G_x$ has a countable paradoxical decomposition.?
\end{proof}
It should be mentioned that the converse of Theorem \ref{thinfinite} does not  hold. For example if $X$ is infinite and $G\curvearrowright X$ is transitive, then for every $x\in X$ we have  $G/G_x$ is infinite whereas $X/G$ is finite.
\begin{remark}
Suppose that $G\curvearrowright X$ is an action of $G$ on $X$ and the sets $G/G_x$ are all countably paradoxical with respect to a fixed couple of translating sets. Then the action $G\curvearrowright X$ is countably pardoxical.
\end{remark}
\begin{proof}
Let every set $G/G_x$ be countably paradoxical with respect to the sequences $\{g_i\}_{i=1}^\infty$ and $\{h_i\}_{i=1}^\infty$. In other words, for every $x\in X$ there are sequences  $\{A^x_i\}_{i=1}^\infty$ and $\{B^x_i\}_{i=1}^\infty$ of susets of $G/G_x$ such that 
\begin{equation*}
G/G_x=\left( \bigsqcup_{i=1}^\infty A^x_i\right)\bigsqcup\left(\bigsqcup_{i=1}^\infty B^x_i\right)
=\bigsqcup_{i=1}^\infty g_i A^x_i
=\bigsqcup_{i=1}^\infty h_i B^x_i. 
\end{equation*}
Setting  $C^x_i=\phi_x(A_i^x)$,  $D^x_i=\phi_x(B_i^x)$ and applying the functions $\phi_x$  one more time, we have  
\begin{eqnarray*}
X&=&\bigsqcup_{x\in X}[x]=\left( \bigsqcup_{i=1}^\infty \left( \bigsqcup_{x\in X} C^x_i\right)\right)\bigsqcup\left(\bigsqcup_{i=1}^\infty \left( \bigsqcup_{x\in X}D^x_i\right)\right)\\
&=&\bigsqcup_{i=1}^\infty g_i \left(\bigsqcup_{x\in X}A^x_i\right)\\
&=&\bigsqcup_{i=1}^\infty h_i \left(\bigsqcup_{x\in X}B^x_i\right). 
\end{eqnarray*}
This completes the proof.
\end{proof}

%
%
%


\end{document}